\documentclass[12pt]{amsart}
\usepackage{amssymb}

\makeatletter

\newtheorem{lemma}{Lemma}[section]
\newtheorem{theorem}[lemma]{Theorem}
\newtheorem{proposition}[lemma]{Proposition}
\newtheorem{corollary}[lemma]{Corollary}

\theoremstyle{definition}
\newtheorem{defn}[lemma]{Definition}
\newtheorem{obs}[lemma]{Observation}

\makeatother

\begin{document}

\title[A finitenes result for commuting squares]{A finiteness result for commuting squares of matrix algebras}         % Enter your title between curly braces
\author[Remus Nicoara]{Remus Nicoara\\University of California Los Angeles}        % Enter your name between curly braces

\begin{abstract}
	We consider a condition for non-degenerate commuting squares of matrix algebras (finite dimensional von Neumann algebras) called the \emph{span condition}, which in the case of the $n$-dimensional standard spin models is shown to be satisfied if and only if $n$ is prime. We prove that the commuting squares satisfying the span condition are isolated among all commuting squares (modulo isomorphisms).  In particular, they are finiteley many for any fixed dimension. Also, we give a conceptual proof of previous constructions of certain one-parameter families of biunitaries.
\end{abstract}
\maketitle

\section{Introduction}

 In this paper we prove some finiteness results for commuting squares of matrix algebras, i.e. finite dimensional von Neumann algebras. Commuting squares were introduced in \cite{Po1}, as invariants and construction data in Jones' theory of subfactors. They encode the generalized symmetries of the subfactor, in a lot of situations being complete invariants (\cite{Po1},\cite{Po2}). In particular, all finite groups and finite dimensional $C^*$-Hopf algebras can be encoded in commuting squares.

One of the simplest examples of commuting squares is
$$\mathfrak{C}=\left(\begin{matrix}
D & \subset{} & M_n(\mathbb{C}) \cr
\cup   &           & \cup \cr
\mathbb{C} & \subset{} & U^*DU
\end{matrix},\tau\right)$$
where $D$ is the algebra of diagonal matrices, $U=(\frac{1}{\sqrt{n}}\epsilon^{(i-1)(j-1)})_{i,j}$ with $\epsilon=cos\frac{2\pi}{n}+isin\frac{2\pi}{n}$, so $U^*DU$ is the algebra of circulant permutation matrices(\cite{Po3}).
 
We call $U$ the \emph{standard biunitary} of order $n$. More generally, one can ask for what unitaries $U$ is $\mathfrak{C}$ a commuting square. The commuting square condition asks that $D,U^*DU$ be orthogonal modulo $\mathbb{C}$, which is equivalent to $U$ having all entries of the same absolute value $1/\sqrt{n}$. Such a matrix is called a \emph{biunitary matrix} or \emph{complex Hadamard matrix}. 

	In \cite{Petrescu} Petrescu showed that, for $n$ positive integer, the standard biunitary of order $n$ is isolated among all normalized biunitary matrices of order $n$ if and only if $n$ is prime. 

We introduce a condition for arbitrary non-degenerate commuting squares, which we call \emph{the span condition}, and prove that it is sufficient to ensure isolation. We show that when the commuting square is given by the standard biunitary of order $n$ the span condition is satisfied if and only if $n$ is prime. Thus our result generalizes Petrescu's finiteness theorem and the span condition can be regarded as a \emph{primeness condition}.

We also show how one can use our theorem to check if a given biunitary is isolated. As an application we show that all circulant biunitaries of order 7 (computed in \cite{Haagerup}) are isolated among all biunitary matrices. 

Conversely, we find sufficient conditions for the span condition to fail and prove that if these conditions are satisfied then there exists a continuum of non-isomorphic commuting squares. 

It is not known if for every $n>5$ prime there exists a one-parameter family of (different) normalized biunitary matrices. For $n=7,13,19,31$ Petrescu found such examples, using a computer; we find a conceptual explanation for these examples. A main point of interest of Petrescu's result is that it might produce examples of one-parameter families of non-isomorphic subfactors of same index $n$ and same graph $A_\infty$.

\section{Preliminaries and a technical result}

We recall the following definition from \cite{Po2}(see also \cite{Po3},\cite{Po1}):
\begin{defn}	A \emph{commuting square} of matrix algebras is a square of inclusions: 
$$
\left(\begin{matrix}
P_{-1} & \subset{} & P_0 \cr
\cup   &           & \cup \cr
Q_{-1} & \subset{} & Q_0 \cr 
\end{matrix},\tau \right)
 $$ 
with $P_0,P_{-1},Q_0,Q_{-1}$ finite dimensional von Neumann algebras (i.e. algebras of the form $\oplus_i \mathbb{M}_{n_i}(\mathbb{C})$, or equivalently *-subalgebras of $\mathbb{M}_n(\mathbb{C}$) for some $n\geq 1$) and $\tau$ a faithful positive trace on $P_0$, $\tau(1)=1$, satisfying the condition:
\begin{equation}\label{csc} E_{P_{-1}}E_{Q_0}=E_{Q_{-1}}\end{equation}
where $E_A=E_A^{P_0}$ denotes the $\tau$-invariant conditional expectation of $P_0$ onto the subalgebra $A\subset P_0$. We say that the commuting square is \emph{non-degenerate} if $P_0=spanP_{-1}Q_0$.
\end{defn}

The following definition is from \cite{Christensen}:

\begin{defn}Let $A$ be a finite dimensional von Neumann algebra with identity $I$ and normalized trace $\tau$. Denote $\mathcal{S}(A)=$ the set of all *-subalgebras of $A$ containing I. For $B_1,B_2\in \mathcal{S}(A)$ and $\delta >0$ we say that $B_1$ is $\delta$-contained in $B_2$ if for every element $x\in B_1$ of $\Vert x \Vert=1$ there exists $y\in B_2$ such that $\Vert x-y \Vert_2 <\delta$. If $B_1$ is $\delta$-contained in $B_2$ and $B_2$ is $\delta$-contained in $B_1$ we write $\Vert B_1-B_2 \Vert_{2,A}<\delta$.
\end{defn}

\begin{obs}\label{isoobs}
Arguments from \cite{Christensen} show that there exists a continous increasing function $f:[0,\infty)\rightarrow [0,\infty), f(0)=0,$ such that if $\delta$ is small and $\Vert B_1-B_2 \Vert_{2,A}<\delta$, then $B_2=Ad(U)(B_1)$ for some unitary element $U\in A$, $\Vert U-I \Vert_2<f(\delta)$ (where $\Vert x \Vert_2 =\tau(x^*x)^{1/2}$). 
\end{obs}

\begin{defn} We say that the commuting square $$\mathfrak{C}=\left(\begin{matrix}
P_{-1} & \subset  & P_0 \cr 
\cup & & \cup \cr 
Q_{-1} & \subset & Q_0 \cr 
\end{matrix},\tau\right)
$$is \textbf{isomorphic} to the commuting square$$\tilde{\mathfrak{C}}=\left(\begin{matrix}
\tilde{P}_{-1} & \subset  & \tilde{P_0} \cr 
\cup & & \cup \cr 
\tilde{Q}_{-1} & \subset & \tilde{Q_0} \cr 
\end{matrix},\tilde{\tau}\right)
$$with trace $\tilde{\tau}$, if there exists a trace-invariant *-isomorphism $\phi:P_0\rightarrow \tilde{P_0}$ such that $\phi(P_{-1})=\tilde{P}_{-1}$, $\phi(Q_{-1})=\tilde{Q}_{-1}$, $\phi(Q_0)=\tilde{Q_0}$
\end{defn}

We can now give the following:

\begin{defn}\label{isodef}We say that the commuting square of matrix algebras $$\mathfrak{C}=\left(\begin{matrix}
P_{-1} & \subset  & P_0 \cr 
\cup & & \cup \cr 
Q_{-1} & \subset & Q_0 \cr 
\end{matrix},\tau\right)
$$is \emph{isolated} if there exists $\delta >0$ such that if 
$$\tilde{\mathfrak{C}}=\left(\begin{matrix}
\tilde{P}_{-1} & \subset  & \tilde{P_0} \cr 
\cup & & \cup \cr 
\tilde{Q}_{-1} & \subset & \tilde{Q_0} \cr 
\end{matrix},\tilde{\tau}\right)
$$ is a commuting square and $\phi:P_0\rightarrow \tilde{P}_0$ a trace-invariant *-isomorphism satisfying $$\Vert \phi(P_{-1})-\tilde{P}_{-1}\Vert_{2,\tilde{P_0}}<\delta, \Vert \phi(Q_{-1})-\tilde{Q}_{-1}\Vert_{2,\tilde{P_0}}<\delta, \Vert \phi(Q_0)-\tilde{Q_0}\Vert_{2,\tilde{P_0}}<\delta$$ then $\tilde{\mathfrak{C}}$ is isomorphic to $\mathfrak{C}$.
\end{defn} 

\noindent For algebras $B\subset A$ we will use the notation: $$B'\cap A=\{a\in A \text{ such that } ab=ba, \forall b\in B\}$$

\begin{lemma}\label{modify} Let $P_0,P_{-1},Q_0,Q_{-1}$ be finite dimensional von Neumann algebras, and $U$ a unitary element of $P_0$ such that
$$
\mathfrak{C}(U)=\left(\begin{matrix}
P_{-1} & \subset  & P_0 \cr 
\cup & & \cup \cr 
Q_{-1} & \subset & U^*Q_0U \cr 
\end{matrix},\tau \right)
$$ 
is a commuting square. Let $q\in Q_0,q'\in Q_0^{'}\cap{}P_{-1},p\in Q_{-1}'\cap P_{-1},p'\in P_{-1}^{'}\cap{}P_0$ be unitary elements. Then $\mathfrak{C}(qq^{'}Upp^{'})$ is a commuting square isomorphic to $\mathfrak{C}(U)$.

\end{lemma}

\begin{proof} Modifying $U$ to the left by $q,q'$ does not change the algebra $U^*Q_0U$ and thus does not change the commuting square: $\mathfrak{C}(qq^{'}Upp^{'})=\mathfrak{C}(Upp^{'})$. By applying $Ad(pp')$ to $\mathfrak{C}(Upp^{'})$ (which leaves $P_0,P_{-1},Q_{-1}$ invariant) we see that $\mathfrak{C}(Upp^{'})$ is isomorphic to $\mathfrak{C}$.

\end{proof}

To check in practical situations if a certain commuting square is isolated, we need the following lemma:

\begin{lemma}\label{isocheck} Let $$\mathfrak{C}=\left(\begin{matrix}
P_{-1} & \subset  & P_0 \cr 
\cup & & \cup \cr 
Q_{-1} & \subset & Q_0 \cr 
\end{matrix},\tau\right)
$$
be a commuting square of finite dimensional von Neumann algebras, with trace $\tau$. $\mathfrak{C}$ is isolated if and only if there exists $\varepsilon>0$ such that if $U\in Q_{-1}'\cap P_0$ is a unitary, $\Vert U-I\Vert_2<\varepsilon$, and
$$\mathfrak{C}(U)=\left(\begin{matrix}
P_{-1} & \subset  & P_0 \cr 
\cup & & \cup \cr 
Q_{-1} & \subset & U^*Q_0U \cr 
\end{matrix},\tau\right)
$$
is a commuting square, then $\mathfrak{C}(U)$ is isomorphic to $\mathfrak{C}$.

\end{lemma}

\begin{proof} We only have to show the implication from right to left. Assume $\mathfrak{C}$ is isolated among commuting squares of the form $\mathfrak{C}(U),U\in Q_{-1}'\cap P_0, \Vert U-I\Vert_2<\varepsilon$ and let $\delta>0$ be such that $f(2f(\delta))+f(\delta)+f(f(\delta)+f(2f(\delta)))<\varepsilon$, where $f$ is as in Observation \ref{isoobs}. We show that $\mathfrak{C},\delta$ satisfy the definition of isolation \ref{isodef}. Assume $\mathfrak{\tilde{C}}$ is $\delta$-close to $\mathfrak{C}$ as in \ref{isodef}. For $\delta$ small the inclusions $\tilde{Q}_{-1}\subset \tilde{P}_{-1}\subset \tilde{P}_0$ and $\phi(Q_{-1})\subset \phi(P_{-1})\subset \tilde{P}_0$ are unitary conjugate. Because our definition of isolation is invariant to isomorphisms of commuting squares, it follows that to check if $\mathfrak{C}$ is isolated it is enough to check isolation among commuting squares of the form:
$$
\begin{matrix}
P_{-1} & \subset  & P_0 \cr 
\cup & & \cup \cr 
Q_{-1} & \subset & \tilde{Q}_0 \cr 
\end{matrix}
$$ 
If $\Vert \tilde{Q}_0-Q_0\Vert_{2,P_0}<\delta$ then using Observation \ref{isoobs} we get $\tilde{Q}_0=U^*Q_0U$ for some unitary $U\in P_0$, $\Vert U-I \Vert_2 <f(\delta)$.

Since $UQ_{-1}U^*\subset Q_0$ and $\Vert UQ_{-1}U^*-Q_{-1} \Vert_{2,Q_0}<2f(\delta)$, Observation \ref{isoobs} implies the existence of a unitary $r_1\in Q_0$, $\Vert r_1-I\Vert_2<f(2f(\delta))$, such that $UQ_{-1}U^*=r_1Q_{-1}r_1^*$. So $Ad(r_1^*U)$ is an isomorphism of $Q_{-1}$ $f(\delta)+f(2f(\delta))$-close to identity, therefore: $Ad(r_1^*U)_{|Q_{-1}}=Ad(r_2)$, for some $r_2\in Q_{-1}$, $\Vert r_2-I\Vert_2 <f(f(\delta)+f(2f(\delta)))$. Thus, by changing $U$ to $r_1^*Ur_2^*$ (which does not change the isomorphism class of the commuting square), we may assume that $U\in{}Q_{-1}^{'}\cap{P_0}$, and since we chosed $\varepsilon>f(2f(\delta))+f(\delta)+f(f(\delta)+f(2f(\delta)))$ we get $\mathfrak{\tilde{C}}$ isomorphic to $\mathfrak{C}$.
 
\end{proof}

According to Lemma \ref{isocheck}, if a commuting square $\mathfrak{C}$ is not isolated then there exists a sequence of unitaries $U_n\rightarrow I$ such that $\mathfrak{C}(U_n)$ are non-isomorphic to $\mathfrak{C}, \forall n\geq 1$. In our main theorem we prove that commuting squares satisfying a certain \emph{span condition} are isolated. To do this, we contradict isolation by assuming the existence of such $U_n$, then we write the commuting square relations for each $n$ and take the "derivative" of this relations along some "direction of convergence" of $U_n$. We want to give a clear meaning to this notion.

Let $P_0$ be a finite dimensional von Neumann algebra and let $U_n=exp(ih_n), n\geq 1$, with $h_n\in P_0$ hermitian non-zero elements converging to $0$. Because of the compactness of the unit ball in the finite dimensional algebra $P_0$, by eventually passing to a subsequence of $\mathbb{N}$ we may assume that $\frac{h_n}{\Vert h_n \Vert}\rightarrow h\in P_0, \Vert h \Vert =1$ (which we will call the \emph{direction of convergence} of $(U_n)_n$). 

Since $\frac{U_n-I}{i\Vert h_n \Vert} \rightarrow h$ as $n\rightarrow\infty$, it follows $\frac{\Vert U_n-I \Vert}{\Vert h_n \Vert} \rightarrow \Vert h \Vert=1$ so:

\begin{equation}\label{h}
h=lim_{n\rightarrow \infty} \frac{U_n-I}{i\Vert U_n-I \Vert}
\end{equation}

The following technical lemma is essential for the proof of the main theorem, giving a way to normalize $h$ by modifying $U_n$ as in Lemma \ref{modify}:

\begin{lemma}\label{norml} With the notations of \ref{modify}, assume that $U_n\in Q_{-1}'\cap P_0, n\geq 1$ are unitary elements converging to $I$ such that $\mathfrak{C}(U_n)$ are commuting squares non-isomorphic to $\mathfrak{C}$.

Then, after replacing $(U_n)_n$ with one of its subsequences, there exist unitaries $q_n\in Q_{-1}'\cap Q_0,q'_n\in Q_0^{'}\cap{}P_{-1},p_n\in Q_{-1}'\cap P_{-1},p'_n\in P_{-1}^{'}\cap{}P_0$ such that:

$$\tilde{U_n}=q_nq_n{'}U_np_n'p_n\rightarrow I, lim_{n\rightarrow \infty} \frac{\tilde{U_n}-I}{i\Vert \tilde{U_n}-I \Vert}=\tilde{h}\in P_0$$ and:

 $$E_{P_{-1}'\cap P_0}(\tilde{h})=E_{Q_0'\cap P_0}(\tilde{h})=E_{Q_{-1}'\cap P_{-1}}(\tilde{h})=E_{Q_{-1}'\cap Q_0}(\tilde{h})=0, [\tilde{h},Q_{-1}]=0$$ 

\end{lemma}

\begin{proof} Let $\mathfrak{X}=\mathfrak{U}(Q_{-1}'\cap Q_0)$ x $\mathfrak{U}(Q_{0}'\cap P_0)$ x $\mathfrak{U}(P_{-1}'\cap P_0)$ x $\mathfrak{U}(Q_{-1}'\cap P_{-1})$ be the set of quadruples of unitaries in the four algebras. $\mathfrak{X}$ being compact in $\Vert$ $ \Vert_2$, for every $n$  there are elements $q_n\in Q_{-1}'\cap Q_0,p_n\in Q_0^{'}\cap{}P_{-1},q'_n\in Q_{-1}'\cap P_{-1},p'_n\in P_{-1}^{'}\cap{}P_0$ that realize the minimum:

$$\Vert q_nq_n{'}U_np_n'p_n-I \Vert_2=\inf_{(q,q',p,p')\in \mathfrak{X}} \Vert qq'U_np'p-I \Vert_2$$

Define $\tilde{U_n}=q_nq_n{'}U_np_n'p_n\rightarrow I$, as for $p=p'=q=q'=I$ we get: $\Vert \tilde{U_n}-I \Vert_2 \leq \Vert U_n-I \Vert_2$. Note that $U_n\not=I$ because the commuting squares were assumed non-isomorphic.Since for every unitary $U$ we have: $\Vert U-I \Vert_2^2=2-2\Re\tau(U)$ (where $\Re\tau$ is the real part of $\tau$), it follows:

$$\Re\tau(\tilde{U_n})\geq \Re\tau(qq'U_np'p), \forall (q,q',p,p')\in \mathfrak{X}$$

Let $\lambda$ be a real number, let $q_0\in Q_{-1}'\cap Q_0$ be a hermitian element, and let $q=exp(i\lambda q_0)q_n$, $q'=q'_n$, $p=p_n$, $p'=p'_n$. Then:

$$\Re\tau(\tilde{U_n})\geq \Re\tau(exp(i\lambda q_0)\tilde{U_n})\Longrightarrow \Re\tau((exp(i\lambda q_0)-I)\tilde{U_n}) \leq 0$$

By dividing with $\lambda > 0$ and taking limit as $\lambda$ approaches 0, we get:
$\Re\tau(iq_0\tilde{U_n})\leq 0$; doing the same but with $\lambda < 0$ we get $\Re\tau(iq_0\tilde{U_n})\geq 0$ so it follows:
$$\Re\tau(iq_0\tilde{U_n})=0$$

Since for hermitians $q_0$ we have $\Re\tau(iq_0I)=0$, we can rewrite the previous equality as:

$$\Re\tau(iq_0(\tilde{U_n}-I))=0$$

Let now (after passing to a subsequence) $\tilde{h}=lim_{n\rightarrow \infty} \frac{\tilde{U_n}-I}{i\Vert \tilde{U_n}-I \Vert} $. Dividing the previous equality by the real number $\Vert \tilde{U_n}-I \Vert$ and taking the limit we get: $$\Re\tau(iq_0(i\tilde{h}))=0 \Longrightarrow \tau(q_0\tilde{h})=0$$
were we used that $\tau(q_0\tilde{h})$ is a real number (since $q_0,\tilde{h}$ are hermitians). Since $Q_{-1}'\cap Q_0$ is the span of its self-adjoint elements, it follows that $E_{Q_{-1}'\cap Q_0}(h)=0$.

Similarly we show that all four expectations are zero (for instance choose $q'$ to be $exp(i\lambda q'_0)q'_n$ and do the same trick, using the fact that $exp(i\lambda q'_0)$ commutes with $q_n$ so it can be moved to the left of the formula for $\tilde{U_n}$).

Since we only modified $U_n$ with elements commuting with $Q_{-1}$, we have $[\tilde{h},Q_{-1}]=0$.

\end{proof}

\section{Span Condition and the main result}

We introduce the \emph{span condition} and show that a commuting square satisfying it is isolated among all commuting square (modulo isomorphisms). In the next lemmas we will often use the following relation that holds true for every $a,b,c\in P_0$: \begin{equation}\label{tcom} \tau([a,b]c)=\tau(a[b,c])=\tau([c,a]b) \end{equation}
as it can be easily checked: $\tau([a,b]c)=\tau(abc-bac)=\tau(abc)-\tau(bac)=\tau(abc)-\tau(acb)=\tau(a[b,c])=\tau(cab)-\tau(acb)$

Before defining the span condition, we present a lemma that justifies it. For $V,W$ vector subspaces of the algebra $P_0$, denote $$V+W=\{v+w| v\in V,w\in W\}$$ $$[V,W]=span\{vw-wv|v\in V,w\in W\}$$
\begin{lemma}\label{span} Let $$\left(\begin{matrix}
P_{-1} & \subset  & P_0 \cr 
\cup & & \cup \cr 
Q_{-1} & \subset & Q_0 \cr 
\end{matrix},\tau\right)$$
be a commuting square with normalized trace $\tau$. Then the vector space $Q_{-1}'\cap P_{-1}+Q_{-1}'\cap Q_0+P_{-1}'\cap P_0+Q_0'\cap P_0$ is orthogonal on $[P_{-1},Q_0]$, with respect to the inner product defined by $\tau$ on $P_0$.

\end{lemma}

\begin{proof} 
Let $p\in P_{-1}$ and $q\in Q_0$. The commuting square condition $$E_{P_{-1}}E_{Q_0}=E_{Q_{-1}}$$ implies $E_{P_{-1}}(q)=E_{Q_{-1}}(q)$ so $E_{P_{-1}}(q-E_{Q_{-1}}(q))=0$, which implies $\tau((q-E_{Q_{-1}}(q))p)=0 $, wich in turn implies
$$\tau(qp)=\tau(E_{Q_{-1}}(q)p)=\tau(E_{Q_{-1}}(q)E_{Q_{-1}}(p))=\tau(qE_{Q_{-1}}(p))$$

\noindent Let $[p_0,q_0]\in [P_{-1},Q_0]$, and $p_1\in Q_{-1}'\cap P_{-1}, q_1\in Q_{-1}'\cap P_{-1}, p'_1\in P_{-1}'\cap P_0, q'_1\in Q_0'\cap P_0$. Using (\ref{tcom}) and $[p_1,p_0]\in P_{-1}$ we obtain:
$$\begin{aligned}\tau([[p_0,q_0]p_1])& =\tau([p_1,p_0]q_0)\\
& =\tau([p_1,p_0]E_{Q_{-1}}(q_0))\\
& =\tau([E_{Q_{-1}}(q_0),p_1]p_0)\\
& =0\end{aligned}$$
since $[E_{Q_{-1}}(q_0),p_1]=0$. Similarly $\tau([p_0,q_0]q_1)=0$. We also have:
$$\tau([p_0,q_0]p_1')=\tau([p_1',p_0]q_0)=0, \tau([p_0,q_0]q_1')=\tau(p_0[q_0,q_1'])=0$$ 
which ends the proof of the lemma.

\end{proof}

\begin{defn}\label{spancond}We say that the commuting square from Lemma \ref{span} satisfies \emph{the span condition} if: $$[P_{-1},Q_0]+(Q_{-1}'\cap P_{-1})+(Q_{-1}'\cap Q_0)+(P_{-1}'\cap P_0)+(Q_0'\cap P_0)=P_0$$

\end{defn}

\begin{obs} Lemma \ref{span} implies that $$dim[P_{-1},Q_0] \leq dim(P_0)-dim(Q_{-1}'\cap P_{-1}+Q_{-1}'\cap Q_0+P_{-1}'\cap P_0+Q_0'\cap P_0)$$ so in some sense the span condition asks for the dimension of the commutator $[P_{-1},Q_0]$ to be maximal.
\end{obs}

The span condition is a reasonable restriction as long as we assume that the commuting square satisfies some non-degeneracy properties, like $dim(P_{-1}'\cap Q_0)=dim(Q_{-1})$, $P_0=spanP_{-1}Q_0$. Indeed, the dimension of $[P_{-1},Q_0]$ is typically big, $P_{-1},Q_0$ are mutually orthogonal (modulo their intersection $Q_{-1}$) and in most of the examples (like the commuting squares associated to groups, Hopf algebras \cite{Szy}, or those corresponding to biunitaries \cite{Petrescu}) their commutants also sit orthogonally (mainly because of the existence of some modular involutions).

We can now prove our main result, which shows that the span condition is sufficient for isolation.

\begin{theorem}\label{spmain}
If the commuting square of finite dimensional von Neumann algebras 
 $$\mathfrak{C}=\left(\begin{matrix}
P_{-1} & \subset  & P_0 \cr 
\cup & & \cup \cr 
Q_{-1} & \subset & Q_0 \cr 
\end{matrix},\tau \right)$$
satisfies the span condition of \ref{spancond}, then $\mathfrak{C}$ is isolated..
\end{theorem}

\begin{proof} Assume, by contradiction, that $\mathfrak{C}$ satisfies the span condition but it is not isolated. According to Lemma \ref{isocheck}, this implies the existence of unitaries $U_n\in P_0,n\geq 1$ converging to $I$ such that:
$$\mathfrak{C}_n=\left(\begin{matrix}
P_{-1} & \subset  & P_0 \cr 
\cup & & \cup \cr 
Q_{-1} & \subset & U_n^*Q_0U_n \cr 
\end{matrix},\tau \right)$$
are commuting squares non-isomorphic to $\mathfrak{C}$. Using Lemma \ref{norml} we may assume:
$$\lim_{n\rightarrow \infty} \frac{U_n-I}{i\Vert U_n-I \Vert}=h\in Q_{-1}'\cap P_0$$
$$E_{P_{-1}'\cap P_0}(h)=E_{Q_0'\cap P_0}(h)=E_{Q_{-1}'\cap P_{-1}}(h)=E_{Q_{-1}'\cap Q_0}(h)=0$$ 

Also:  $$\lim_{n\rightarrow \infty} \frac{U_n^*-I}{i\Vert U_n-I \Vert}=-h$$

Let $p\in P_{-1}$ such that $E_{Q_{-1}}(p)=0$ and let $q\in Q_0$. The commuting square condition implies $E_{U_n^*Q_0U_n}(p)=E_{Q_{-1}}(p)=0$, so $$\tau(pU_n^*qU_n)=0=\tau(pq) \Longrightarrow \tau(p(U_n-I)^*qU_n)+\tau(pq(U_n-I))=0$$
Dividing by $i\Vert U_n-I \Vert$ and taking the limit as $n\rightarrow \infty$ it follows
$$\tau(p(-h)q)+\tau(pqh)=0 \Longrightarrow \tau([p,q]h)=0$$
Thus, the $h$ is orthogonal on all vectors $[p,q]$ with $E_{Q_{-1}}(p)=0$. We show that $h$ is in fact orthogonal on all vectors in $[P_{-1},Q_0]$. Indeed, if $p_1$ is an arbitrary element of $P_{-1}$, using $E_{Q_{-1}}(p_1-E_{Q_{-1}}(p_1))=0$ we get:
$$\begin{aligned}\tau([p_1,q]h) & =\tau([p_1-E_{Q_{-1}}(p_1),q]h+[E_{Q_{-1}}(p_1),q]h)\\
& =0+\tau([E_{Q_{-1}}(p_1),q]h)\\
& =\tau([h,E_{Q_{-1}}(p_1)]q)\\
& =0\end{aligned}$$
We used formula (\ref{tcom}) and $[h,Q_{-1}]=0$. This shows that $h$ is orthogonal on $[P_{-1},Q_0]$. As $h$ is also orthogonal on the algebras $Q_{-1}'\cap P_{-1},Q_{-1}'\cap Q_0,P_{-1}'\cap P_0,Q_0'\cap P_0 $ it follows that if the span condition holds we must have $E_{P_0}(h)=0$ so $h=0$, which contradicts $\Vert h \Vert =1$, contradiction that ends the proof

\end{proof}

\begin{corollary} For every $N\geq 2$ there are only finitely many (up to isomorphisms) commuting squares $\mathfrak{C}$ with $dim(P_0)=N$, satisfying the span condition.

\end{corollary}
$$ $$

\section{Existence of one-parameter families of non-isomorphic commuting squares}

In the previous section we've showed that the span condition is sufficient for isolation, but we did not discuss wether it is also necessary. We give partial converses to Theorem \ref{spmain}, which consider some of the simplest cases in which the span condition fails. The next theorem shows that one can construct a continuum of commuting squares if there exist two non-trivial elements $p_0\in P_{-1},q_0\in Q_0$ that commute.

\begin{theorem}\label{constr1} Let  $$\mathfrak{C}=\left(\begin{matrix}
P_{-1} & \subset  & P_0 \cr 
\cup & & \cup \cr 
Q_{-1} & \subset & Q_0 \cr 
\end{matrix},\tau \right)$$
be a commuting square of finite dimensional von Neumann algebras, and assume there exist hermitian elements $p_0\in Q_{-1}'\cap P_{-1},q_0\in Q_{-1}'\cap Q_0$, that are not in $Q_{-1}$, such that $p_0q_0-q_0p_0=0$. If $U_t=exp(itp_0q_0),t\in\mathbb{R}$, then
$$\mathfrak{C}_t=\left(\begin{matrix}
P_{-1} & \subset  & P_0 \cr 
\cup & & \cup \cr 
Q_{-1} & \subset & U_t^*Q_0U_t \cr 
\end{matrix},\tau \right)$$
is a one-parameter family of commuting squares.
\end{theorem}

\begin{proof} We show that the commuting square condition holds for each $t$. Let $p\in P_{-1}$ such that $E_{Q_{-1}}(p)=0$, and $q\in Q_0$. We need to show that $\tau(pU_t^*qU_t)=0$. Writing $U_t=exp(ipqt)=\sum_{k}\frac{i^kt^k}{k!}p^kq^k$ we have: 

$$\begin{aligned}
\tau(pU_t^*qU_t) & =\sum_{k,l} \frac{(-1)^li^{k+l}t^{k+l}}{k!l!}\tau(pp_0^lq_0^lqq_0^kp_0^k)\\
& = \sum_{k,l} \frac{(-1)^li^{k+l}t^{k+l}}{(k+l)!}C_{k+l}^l\tau(p_0^kpp_0^lq_0^lqq_0^k)\\
& = \sum_{n}\sum_{k+l=n} \frac{(-1)^li^{n}t^{n}}{n!}C_{n}^l\tau(p_0^kpp_0^lq_0^lqq_0^k)\\
& =\sum_{n} (\frac{i^{n}t^{n}}{n!}\tau(p_0^npq_0^nq)(\sum_{k+l=n}(-1)^lC_{n}^l))\\
& =\tau(pq)=\tau(p)\tau(q)\\
& =0 \\
\end{aligned}$$
We used: $$\begin{aligned}\tau(p_0^kpp_0^lq_0^lqq_0^k)& =\tau(E_{Q_{-1}}(p_0^kpp_0^l)E_{Q_{-1}}(q_0^lqq_0^k))\\
& =\tau(E_{Q_{-1}}(p_0^np)E_{Q_{-1}}(q_0^nq))\\
& =\tau(p_0^npq_0^nq) \end{aligned}$$
since $p_0\in Q_{-1}'\cap P_{-1},q_0\in Q_{-1}'\cap Q_0$. We also used: $\sum_{l}(-1)^lC_{n}^l=0$ for $n \geq 1$. 

\end{proof}

If $p,q$ are projections then the unitaries in Theorem \ref{constr1} can be written as $U(\lambda)=I+(\lambda-1)pq$, $\lambda=e^{it}\in\mathbb{T}$. This justifies the class of unitaries we construct in the next theorem, that aplies to situations when there exists a linear dependence relation between 2 commutators in the span.

\begin{theorem}\label{constr2}  Let  $$\mathfrak{C}=\left(\begin{matrix}
P_{-1} & \subset  & P_0 \cr 
\cup & & \cup \cr 
Q_{-1} & \subset & Q_0 \cr 
\end{matrix},\tau \right)$$
be a commuting square of finite dimensional von Neumann algebras, and assume there exist orthogonal projections $p_1,p_2\in Q_{-1}'\cap P_{-1}$ and orthogonal projections $q_1,q_2\in Q_{-1}'\cap Q_0$, that are not in $Q_{-1}$, satisfying $[p_1,q_1]-[p_2,q_2]=0$. Let $$U(\lambda)=I+(\lambda-1)p_1q_1+(\bar{\lambda}-1)p_2q_2$$ 
for $\lambda\in\mathbb{T}$. Then
$$\mathfrak{C}_\lambda=\left(\begin{matrix}
P_{-1} & \subset  & P_0 \cr 
\cup & & \cup \cr 
Q_{-1} & \subset & U(\lambda)^*Q_0U(\lambda) \cr 
\end{matrix},\tau \right)$$
is a one-parameter family of commuting squares.
\end{theorem}

\begin{proof}

Since $$[p_1,q_1]-[p_2,q_2]=0 \Longrightarrow p_1q_1+q_2p_2=p_2q_2+q_1p_1$$ multiplying with $p_1$ to the right, and then with $p_2$ to the left, we have:
$$p_1q_1p_1=p_2q_2p_1+q_1p_1, p_2q_2p_1+p_2q_1p_1=0$$
Similary, by multiplying with $p_2$ to the left we have:
$$p_2q_2p_2=p_2q_2+p_2q_1p_1$$
and summing up the last relations:
\begin{equation}\label{p1q1p1}\begin{aligned}p_1q_1p_1+p_2q_2p_2 &=p_2q_2p_1+q_1p_1+p_2q_2+p_2q_1p_1\\
&=q_1p_1+p_2q_2+(p_2q_2p_1+p_2q_1p_1)\\
&=q_1p_1+p_2q_2\end{aligned}\end{equation}

We now show that $U(\lambda)$ is a unitary:
$$\begin{aligned}U(\lambda)U(\lambda)^*&=(I+(\lambda-1)p_1q_1+(\bar{\lambda}-1)p_2q_2)(I+(\bar{\lambda}-1)q_1p_1+(\lambda-1)q_2p_2)\\
&=I+(\lambda-1)(p_1q_1+q_2p_2)+(\bar{\lambda}-1)(p_2q_2+q_1p_1)\\
&\quad+(\lambda-1)(\bar{\lambda}-1)(p_1q_1p_1+p_2q_2p_2)\\
&=I\end{aligned}$$

We used: $q_1p_1+p_2q_2=p_2q_2+q_1p_1,(\lambda-1)(\bar{\lambda}-1)=-(\lambda-1)-(\bar{\lambda}-1)$ and equation (\ref{p1q1p1}):
$$p_1q_1p_1+p_2q_2p_2=q_1p_1+p_2q_2$$

Let's now check that $\mathfrak{C}(U)$ is a commuting square: for $p\in P_{-1}$ with $E_{Q_{-1}}(p)=0$ and $q\in Q_0$ we have:

$$\begin{aligned}\tau (pU(\lambda)qU(\lambda)^*)&=\tau(pq)\\
&\quad+(\lambda-1)\tau (pp_1q_1q+pqq_2p_2)+(\bar{\lambda}-1)\tau(pp_2q_2q+pqq_1p_1)\\
&\quad+(\lambda-1)(\bar{\lambda}-1)\tau(pp_1q_1qq_1p_1+pp_2q_2qq_2p_2)\\
&\quad+(\lambda-1)^2\tau(pp_1q_1qq_2p_2)+(\bar{\lambda}-1)^2\tau(pp_2q_2qq_1p_1)\end{aligned}$$

But $\tau(pp_1q_1qq_2p_2)=\tau(p_2pp_1q_1qq_2)=\tau(E_{Q_{-1}}(p_2pp_1)E_{Q_{-1}}(q_1qq_2))=0$, because $E_{Q_{-1}}(q_1qq_2)=E_{Q_{-1}}(qq_2q_1)=0$, since $q_2q_1=0$ and $[q_1,Q_{-1}]=0$. Similarly $\tau(pp_2q_2qq_1p_1)=0$. Also:
$$\begin{aligned}\tau(pp_1q_1q+pqq_2p_2)&=\tau(pp_2q_2q+pqq_1p_1)\\
&=\tau(pp_1q_1qq_1p_1+pp_2q_2qq_2p_2)\\
&=\tau(pq(q_1p_1+q_2p_2))\end{aligned}$$

Indeed: 
$$\begin{aligned}\tau(pp_1q_1q+pqq_2p_2)&=\tau(pp_1q_1q)+\tau(qq_2p_2p)\\
&=\tau(E_{Q_{-1}}(pp_1)E_{Q_{-1}}(q_1q))+\tau(E_{Q_{-1}}(qq_2)E_{Q_{-1}}(p_2p))\\
&=\tau(E_{Q_{-1}}(p_1p)E_{Q_{-1}}(qq_1))+\tau(E_{Q_{-1}}(p_2p)E_{Q_{-1}}(qq_2))\\
&=\tau(p_1pqq_1)+\tau(p_2pqq_2)\\
&=\tau(pq(q_1p_1+q_2p_2))\end{aligned}$$
and the other equalities follow similarly. Thus, using $(\lambda-1)(\bar{\lambda}-1)+(\lambda-1)+(\bar{\lambda}-1)=0$ we have: 
$$\tau (pU(\lambda)qU(\lambda)^*)=0$$
which ends the proof.

\end{proof}

\section{Comments on Petrescu's results}

We discuss consequences of the theorems from the previous sections for commuting squares of the form:
$$
\left(\begin{matrix}
D & \subset{} & \mathbb{M}_n(\mathbb{C}) \cr
\cup   &           & \cup \cr
\mathbb{C} & \subset{} & U^*DU
\end{matrix}, \tau\right)$$  
with $D$=the diagonal matrices, $U$ unitary in $\mathbb{M}_(\mathbb{C})$ and $\tau=\frac{1}{n}Tr$ the normalized trace. 

Denote by $(A_{i,j})_{i,j}$ the matrix units of $M_n(\mathbb{C})$, $A_{i,j}=$the matrix having 1 at the intersection of the $i^{th}$ row and $j^th$ column, and only zeros at the other positions. Also, let $D_k=A_{k,k},k=1,...,n$ be an orthogonal basis of $D$. 

If $U=(u_{ij})_{1\leq i,j\leq n}$ the commuting square condition $\tau(D_iU^*D_jU)=\tau(D_i)\tau(D_j)$ can be rewritten as $\overline{u}_{ji}u_{ij}=1/n$. Thus it amounts to all entries of $U$ having the same absolute value $1/\sqrt{n}$. Such a $U$ is called a \emph{biunitary matrix} or complex Hadamard matrix. We say that two biunitaries are \emph{equivalent} if the corresponding commuting squares are isomorphic. For every $n$ there exists at least one biunitary of order $n$: $U=\frac{1}{\sqrt{n}}(\epsilon^{(i-1)(j-1)})_{i,j}, \epsilon=cos\frac{2\pi}{n}+isin\frac{2\pi}{n}$, called the \emph{standard biunitary} of order $n$.

We can apply Theorem \ref{spmain} to commuting squares given by biunitary matrices. Since the algebras $D$ and $U^*DU$ are abelian and orthogonal modulo their intersection $\mathbb{C}I$ the span condition becomes: 
$$dim([D,U^*DU])=n^2-2n+1$$

Thus, we have the following:
\begin{proposition}\label{isobiunit}If $U\in M_n(\mathbb{C})$ is a biunitary matrix such that $dim([D,U^*DU])=n^2-2n+1$ then $U$ is isolated among all biunitaries (up to equivalence).

\end{proposition}

\begin{corollary}[Petrescu's Theorem] The standard biunitary of order n is isolated iff n is prime.
\end{corollary}
\begin{proof}

Assume n is prime and let $U=(\epsilon^{(i-1)(j-1)})_{i,j}$ be the standard biunitary matrix of order n. $U^*DU=S$ is the algebra of circulant permutation matrices. $S_k=\sum_{i}A_{i,i+k}, k=1..n$ give a basis for $S$ (all the indices are considered modulo $n$).

$X_{k,l}=[D_k,S_l]=A_{k,k+l}-A_{k-l,k}$ is a set of generators for $[D,U^*DU]$. If for some complex numbers $c_{k,l}$ we have: $$\sum_{k,l}c_{k,l}X_{k,l}=0 $$
It follows:$$\sum_{k,l}(c_{k,l}A_{k,k+l}-c_{k,l}A_{k-l,k})=0 \Longrightarrow \sum_{i,j} (c_{i,j-i}-c_{j,j-i})A_{i,j}=0$$
so  $c_{i,j-i}=c_{j,j-i}$ and if we denote by $s=j-i$ we have $c_{i,s}=c_{i+s,s}$ so $c_{i,s}=c_{i+ms,s}, \forall m=0,1,...,n-1$. Since $n$ is prime, for $s$ nonzero the elements $0,s,2s,...,(n-1)s$ cover all possible residues mod $n$, so $c_{i,s}=c_{0,s}$. 

Thus the dimension of the kernel of the linear transformation $$(c_{k,l})_{k,l}\rightarrow \sum_{k,l}c_{k,l}X_{k,l}$$ is $(2n-1)$, so its range has dimension $n^2-2n+1$, which shows that the span condition holds.

Conversely, if $n$ is not prime, $n=n_1n_2$ with $n_1,n_2>1$, then $p=\sum_{j \leq n_2} A_{jn_1,jn_1}\in D$ and $q=\sum_{i,j} A_{j,in_1+j}\in U^*DU$ commute so by Theorem \ref{constr1} we can construct a one-parameter family of biunitaries $U(t)=exp(ipqt)$.
\end{proof}

For $n=5$ the standard biunitary is the only biunitary, as proven by U.Haagerup (\cite{Haagerup}). For all $n>5$ prime there exists at least another biunitary which is a circulant matrix (\cite{Bj},\cite{HJ},\cite{MW}), and for every $n$ non prime one can easily construct infinitely many biunitaries.

 S.Popa conjectured that for every $n>5$ prime there exist only finitely many normalized biunitaries (\cite{Po3}). Surprisingly, this turned out to be false: one-parameter families of normalized biunitaries where constructed by M.Petrescu for $n=7,13,19,31,79$ (\cite{Petrescu}). A main point of interest in this result is that it can produce one-parameter families of non-isomorphic subfactors with the same graph, conjectured to be $A_\infty$. While Petrescu's examples have been constructed using the computer, we give a conceptual proof of their existence as a consequence of Theorem \ref{constr2} . We will work the details for one of the two examples for $n=7$, the other examples having similar proofs.

\begin{corollary}[Petrescu's biunitaries]Let $\lambda\in \mathbb{T}$, $w=cos\frac{2\pi}{6}+isin\frac{2\pi}{6}$ and
\begin{equation}\label{Petex7}U(\lambda)=\frac{1}{\sqrt{7}}\begin{pmatrix}
\lambda w & \lambda w^4 & w^5 & w^3 & w^3 & w & 1 \cr
\lambda w^4 & \lambda w & w^3 & w^5 & w^3 & w & 1 \cr
w^5 & w^3 & \overline{\lambda} w & \overline{\lambda} w^4 & w & w^3 & 1 \cr
w^3 & w^5 & \overline{\lambda} w^4 & \overline{\lambda} w & w & w^3 & 1 \cr
w^3 & w^3 & w & w & w^4 & w^5 & 1 \cr
w & w & w^3 & w^3 & w^5 & w^4 & 1 \cr
1 & 1 & 1 & 1 & 1 & 1 & 1 \cr
\end{pmatrix}
\end{equation} Then $U(\lambda)$ is a 1-parameter family of (non-equivalent) biunitaries.
\end{corollary}

\begin{proof}Let $U=U(1)$, $P_0=\mathbb{M}_n(\mathbb{C}),P_{-1}=D, Q_0=U^*DU, Q_{-1}=\mathbb{C}$ and
 $$p_1=A_{1,1}+A_{2,2}, p_2=A_{3,3}+A_{4,4}\in P_{-1}$$
 $$q_1=U^*(A_{1,1}+A_{2,2})U, q_2=U^*(A_{3,3}+A_{4,4})U \in Q_0$$
It is easy to check that:
 $$[p_1,q_1]-[p_2,q_2]=0, p_1p_2=q_1q_2=0$$
Thus we are in the conditions of Theorem \ref{constr2}, so $$U(\lambda)=(I+(\lambda -1)p_1q_1+(\bar{\lambda}-1)p_2q_2)U$$ are biunitaries for all $\lambda$ complex numbers of absolute value 1. One can easily verify that $U(\lambda)$ are the biunitaries from (\ref{Petex7}).

\end{proof}

\begin{obs} One can try to find more examples of biunitaries using the following algorithm: fix $p_1,p_2,p_3,p_4\in D$, with $p_1p_2=p_3p_4=0$, and find (with the help of a computer) the local minimum of the function: $$U \rightarrow \Vert UU^*-I \Vert + \Vert [p_1,U^*p_3U]+[p_2,U^*p_4U] \Vert$$
by starting with an arbitrary value of $U$ and taking succesive iterations. If the local minimum is $0$ then one can construct from $U$ a one-parameter family of biunitaries as before.
\end{obs}
\section{Application to circulants}

A matrix $S$ is called \emph{circulant} if all rows are obtained from consecutive circular permutations of the first row, i.e. $S=(s_{j-i})_{i,j\in \mathbb{Z}/n\mathbb{Z}}$. The problem of classifying circulant biunitaries is equivalent to Bjorck's problem of classifying \emph{cyclic n-roots} (\cite{Bj}). For every $n$ prime there exists at least one circulant biunitary, defined as: $s_i=1/\sqrt{n}$ for $i$ quadratic residue modulo $n$ and $s_i=a$ in rest, where $a$ is the root of a certain quadratic equation over $\mathbb{Q}$ (\cite{Bj},\cite{HJ},\cite{MW}). 

U.Haagerup proved that for every $n$ prime there exist finitely many circulant biunitaries (up to equivalence). We conjecture that in fact every circulant biunitary satisfies the span condition, and thus is isolated among all biunitary matrices. We show this is true for $n=7$, by using Haagerup's classification of biunitaries of order 7 (\cite{Haagerup}). We give an algorithm that can be used, more generally, to check if a given $U$ satisfies the span condition, and thus is isolated among all normalized biunitary matrices.
\begin{proposition} If $U$ is a circulant biunitary matrix of order $7$ then $U$ is isolated among all biunitary matrices.
\end{proposition}
\begin{proof}Let $U=(u_{ij})_{i,j\in \mathbb{Z}_7}$ be a circulant biunitary of order $7$. According to Proposition \ref{isobiunit} it is enough to show that $dim(span[D,U^*DU])=7^2-27+1=36$. Let $D_i, i=0...6$ be a basis for $D$, where $D_i$ is the diagonal having $1$ on the $i+1$ position on the diagonal and only $0$'s on the other positions. Let: $a_{kl}^{ij}=[D_i,U^*D_jU]_{kl}=(D_iU^*D_jU-U^*D_jUD_i)_{kl}=\delta_i^k \overline{u}_{ji}u_{jl}-\delta^l_i\overline{u}_{jk}u_{ji}=(\delta_i^k-\delta_i^l)\overline{u}_{jk}u_{jl}$.

Let $A$ be the $49\times 49$ matrix given by $A_{(i,j),(k,l)}=a_{kl}^{ij}$. We need to check that $rank(A)=36$. Because $\sum_{i}a_{kl}^{ij}=\sum_{i}[D_i,U^*D_jU]_{kl}=[I,U^*D_jU]_{kl}=0$ and similarly $\sum_{j}a_{kl}^{ij}=0$, to find the rank of $A$ we may remove the $13$ rows of $A$ indexed after $(i,0),(0,j), 0\leq i,j\leq 6$. Similar arguments show that we can remove the $13$ columns of $A$ indexed after $(i,i),(0,j), 0\leq i,j\leq 6$. If we denote by $M$ the $36\times 36$ matrix left, we need to check that $det(M)\not=0$. 

We include this computation for one of the circulant matrices of order 7 (of the type described at the beginning of this section, for $a=-3/4+i\sqrt{7}/4$). Similarly, we checked that all circulant biunitaries of order 7 (computed in \cite{Haagerup}) are isolated among all biunitaries. The pairs $(i,j),0\leq i,j\leq 6$ are encoded in the vector w, and the selection of rows and columns of $A$ is encoded in v.
$$ $$
with(LinearAlgebra):

ID:=Matrix(7,7,shape=identity);a:=-3/4+I*sqrt(7)/4;

U:=1/sqrt(7)Matrix([[1,1,1,a,1,a,a],[a,1,1,1,a,1,a],[a,a,1,1,1,a,1],
[1,a,a,1,1,1,a],[a,1,a,a,1,1,1],[1,a,1,a,a,1,1],[1,1,a,1,a,a,1]]);

A:=(i,j,k,l)-$>$(ID[k,i]-ID[l,i])*conj(U[j,k])*U[j,l];

w:=m-$>$(1+((m-(1+((m-1) mod 6)))/6),1+((m-1) mod 6));

 v := array(1..72,[1,2,1,3,1,4,1,5,1,6,1,7,2,1,2,3,2,4,2,5,2,6,2,7,3,1,3,

2,3,4,3,5,3,6,3,7,4,1,4,2,4,3,4,5,4,6,4,7,5,1,5,2,5,3,5,4,5,6,5,7,6,1,

6,2,6,3,6,4,6,5,6,7]);

 g:=(n,m)-$>$A(w(m),v[2*n-1],v[2*n]);

 M:=Matrix(36,g); Determinant(M);

\end{proof}

\bibliographystyle{amsalpha}

\begin{thebibliography}{10}

\bibitem[Bj]{Bj}G.Bjorck, \emph{Functions of modulus 1 on $\mathbb{Z}_n$, whose Fourier transforms have constant modulus, and "cyclic n-roots"}, Recent Advances in Fourier Analysis and its applications, NATO Adv. Sci. Inst. Ser. C: Math. Phys. Sci., Kluwer \textbf{315} (1990), 131-140.

\bibitem[Chr]{Christensen}E.Christensen, \emph{Subalgebras of a finite algebra}, Mathematische Annalen \textbf{243}, 17-29 (1979)

\bibitem[H]{Haagerup}U.Haagerup, \emph{Orthogonal maximal abelian *-subalgebras of the $n\times n$ matrices and cyclic $n$-roots}, Operator Algebras and Quantum Field Theory (ed. S.Doplicher et al.), International Press (1997),296-322

\bibitem[HJ]{HJ}P.de la Harpe and V.F.R.Jones, \emph{Paires de sous-algebres semi-simples et graphes fortement reguliers}, C.A. Acad. Sci. Paris \textbf{311}, serie I (1990), 147-150

\bibitem[Jon]{Jones}V.F.R.Jones, \emph{Index for subfactors}, Invent. Math \textbf{72} (1983),
  1--25.

\bibitem[JS]{JS}V.F.R.Jones and V.S.Sunder, \emph{Introduction to subfactors}, London Math. Soc. Lecture Notes Series \textbf{234}, Cambridge University Press, 1997.

\bibitem[MW]{MW}A.Munemasa and Y.Watatani, \emph{Orthogonal pairs of *-subalgebras and Association Schemes} C.R. Acad. Sci. Paris \textbf{314}, serie I (1992), 329-331.

\bibitem[Pe]{Petrescu}M.Petrescu, \emph{Existence of continuous families of complex Hadamard matrices of certain prime dimensions and related results}, PhD thesis, University of California Los Angeles, 1997.

\bibitem[Po1]{Po1}S.Popa, \emph{Classification of subfactors : the reduction to commuting squares}, Invent. Math., \textbf{101}(1990),19-43

\bibitem[Po2]{Po2}S.Popa, \emph{Classification of amenable subfactors of type II}, Acta Mathematica, \textbf{172}, 163-255 (1994).

\bibitem[Po3]{Po3}S.Popa, \emph{Othogonal pairs of *-subalgebras in finite von Neumann algebras}, J. Operator Theory \textbf{9}, 253-268 (1983)

\bibitem[Szy]{Szy}W.Szymanski, \emph{Finite index subfactors and Hopf algebra crossed products}, Proc. Amer. Math. Soc., \textbf{120}(1994),519-528.


\end{thebibliography}

\end{document}